\documentclass[envcountsect]{svmult}
\usepackage[T1]{fontenc}
\usepackage[latin9]{inputenc}
\usepackage{prettyref}
\usepackage{url}
\usepackage{amsmath}
\usepackage{amssymb}
\usepackage{esint}
\usepackage[unicode=true,
 bookmarks=true,bookmarksnumbered=false,bookmarksopen=false,
 breaklinks=false,pdfborder={0 0 1},backref=false,colorlinks=false]
 {hyperref}
\usepackage{breakurl}

\makeatletter
\newenvironment{svmultproof2}{\begin{proof}}{\smartqed\qed\end{proof}}

\usepackage{helvet}
\usepackage{courier}

\newcommand{\isdef}{\ensuremath{\stackrel{\text{def}}{=}}}

\newcommand{\SP}{\ensuremath{{\cal S_P}}}
\newcommand{\SPD}{\ensuremath{{\cal S}_{{\cal P},\Delta}}}

\newrefformat{prob}{Problem \ref{#1}}
\newrefformat{prop}{Proposition \ref{#1}}
\newrefformat{def}{Definition \ref{#1}}
\newrefformat{sub}{Section \ref{#1}}
\newrefformat{alg}{Algorithm \ref{#1}}
\newrefformat{con}{Conjecture \ref{#1}}
\newrefformat{cor}{Corollary \ref{#1}}
\newrefformat{rem}{Remark \ref{#1}}
\newrefformat{app}{Appendix \ref{#1}}
\newrefformat{ex}{Example \ref{#1}}
\newrefformat{enu}{Item \ref{#1}}

\DeclareMathOperator{\Op}{\mathfrak{D}}

\usepackage{arydshln}
\usepackage{multirow}

\makeatother

\begin{document}

\title*{Taylor Domination, Difference Equations, and Bautin Ideals}

\author{Dmitry Batenkov and Yosef Yomdin}

\institute{Dmitry Batenkov \at Department of Computer Science, Technion - Israel Institute of Technology, Haifa 32000, Israel.
\email{batenkov@cs.technion.ac.il}\and Yosef Yomdin \at Department
of Mathematics, Weizmann Institute of Science, Rehovot 76100, Israel.
This author is supported by {\small ISF, Grants No. 639/09 and 779/13, and by the
Minerva foundation.} \email{yosef.yomdin@weizmann.ac.il}}

\maketitle

Abstract. {\small We compare three approaches to studying the behavior of an analytic function $f(z)=\sum_{k=0}^\infty a_kz^k$
from its Taylor coefficients. The first is ``Taylor domination'' property for $f(z)$ in the complex disk $D_R$, which is an
inequality of the form
\[
|a_{k}|R^{k}\leq C\ \max_{i=0,\dots,N}\ |a_{i}|R^{i}, \ k \geq N+1.
\]
The second approach is based on a possibility to generate $a_k$ via recurrence relations. Specifically, we consider linear
non-stationary recurrences of the form
\[
a_{k}=\sum_{j=1}^{d}c_{j}(k)\cdot a_{k-j},\ \ k=d,d+1,\dots,
\]
with uniformly bounded coefficients.

In the third approach we assume that $a_k=a_k(\lambda)$ are polynomials in a
finite-dimensional parameter $\lambda \in {\mathbb C}^n.$ We study ``Bautin ideals'' $I_k$ generated by
$a_{1}(\lambda),\ldots,a_{k}(\lambda)$ in the ring ${\mathbb C}[\lambda]$ of polynomials in $\lambda$.

\smallskip

These three approaches turn out to be closely related. We present some results and questions in this direction.
}

\bigskip
\bigskip
\bigskip
\bigskip

\global\long\def\O{\Omega}
\global\long\def\C{\mathbb{C}}
\global\long\def\N{\mathbb{N}}
\global\long\def\R{\mathbb{R}}
\global\long\def\e{\varepsilon}
\global\long\def\Z{\mathbb{Z}}

\section{Introduction}

Let $f(z)$ be an analytic function represented in a disk $D_R$ of radius $R$, centered at the origin by a convergent Taylor
series $f(z)=\sum_{k=0}^{\infty}a_{k}z^{k}$. We assume that the Taylor coefficients $a_k$ are explicitly known, or, at least,
can be recovered by a certain explicit procedure. In many cases this is the only analytic information we possess on $f(z)$.
A notorious example is the ``Poincar\'e first return mapping'' of a non-linear ordinary differential equation (see
Section \ref{Abel} below). We would like to investigate the behavior of $f(z)$ on the base of what we know, i.e. properties of the sequence $a_k$. In particular,
we would like to bound from above the possible number of the zeroes of $f(z)$. In Section \ref{TD.zeroes} below we provide an explicit connection between the two.

The main goal of this paper is to compare three approaches to the above question: the first is ``Taylor domination'', which
is a bound on all the Taylor coefficients $a_k$ of $f$ through the first few of them. The second is a possibility to
generate $a_k$ via Difference Equations, specifically, by linear non-stationary homogeneous recurrence relations of a fixed
length, with uniformly bounded coefficients. In the third approach we assume that $a_k=a_k(\lambda)$ are polynomials in a
finite-dimensional parameter $\lambda,$ and study ``Bautin ideals'' $I_k$ generated by $a_{1}(\lambda),\ldots,a_{k}(\lambda)$
in the ring of polynomials in $\lambda$.

The main facts which we present are the following:

\smallskip

\noindent 1. A sequence $a_0,a_1,\ldots$ can be obtained as a solution of a non-stationary linear homogeneous recurrence
relations of a fixed length, with uniformly bounded coefficients, if and only if each its subsequence $a_m,a_{m+1},\ldots$
possesses an appropriate Taylor domination property.

\smallskip

\noindent 2. A sequence $a_{1}(\lambda), a_{2}(\lambda)\ldots,$ of polynomials in $\lambda$, under some natural assumptions,
possesses a uniform in $\lambda$ Taylor domination, with the parameters determined through the algebra of the Bautin ideals.

\smallskip

\noindent 3. If a sequence $a_{1}(\lambda), a_{2}(\lambda)\ldots,$ of polynomials in $\lambda$ is produced by {\it an algebraic}
recurrence relation, its Bautin ideals can be computed explicitly. We discuss briefly the difficulties which arise for
{\it differential-algebraic} recurrences (like in case of Poincar\'e mapping of Abel differential equation).

\smallskip

In our discussion we present some new specific results, and many known ones (some of them very recent, some pretty old). We
believe that a general picture of interconnections between Taylor domination, recurrence relations, and Bautin ideals, given
in this paper is new and may be instructive in further developments.

\section{Taylor domination and counting zeroes}\label{TD.zeroes}

``Taylor domination'' for an analytic function $f(z)=\sum_{k=0}^{\infty}a_{k}z^{k}$ is an explicit bound of all its
Taylor coefficients $a_k$ through the first few of them. This property was classically studied, in particular, in relation
with the Bieberbach conjecture, which asserts that for univalent $f$ it always holds that $|a_k| \leq k|a_1|$. See
\cite{bieberbach1955analytische,biernacki1936fonctions,hayman1994multivalent} and references therein. To give an
accurate definition, let us assume the radius of convergence of the Taylor series for $f$ is $\hat R$, \
$0<\hat{R}\leqslant+\infty$.

\begin{definition}
\label{def:domination}Let a positive \emph{finite} $R\leq\hat{R},$
a natural $N$, and a positive sequence $S\left(k\right)$ of a subexponential growth be fixed. The function $f$ is said to
possess an $(N,R,S(k))$ - Taylor domination property if for each $k\geq N+1$ we have
\[
|a_{k}|R^{k}\leqslant S(k)\ \max_{i=0,\dots,N}|a_{i}|R^{i}.
\]
If $S(k)\equiv C$ we call this property $(N,R,C)$-Taylor domination.
\end{definition}
The parameters $(N,R,S(k))$ of Taylor domination are not defined uniquely. In fact, the following easy result of \cite{Bat.Yom.2014}
shows that each nonzero analytic function $f$ possesses this property:

\begin{proposition}\label{prop:Root}[Proposition 1.1, \cite{Bat.Yom.2014}]
If $0<\hat{R}\leqslant+\infty$
is the radius of convergence of $f\left(z\right)=\sum_{k=0}^{\infty}a_{k}z^{k}$, with $f\not\not\equiv0$, then for each
finite and positive $0<R\leqslant\hat{R}$, $f$ satisfies the $\left(N,R,S\left(k\right)\right)$-Taylor domination
property with $N$ being the index of its first nonzero Taylor coefficient, and $S\left(k\right)= R^k|a_k|(|a_{N}|R^{N})^{-1},$
for $k > N$.
\end{proposition}
Consequently, the Taylor domination property becomes really interesting only for those {\it families} of analytic functions
$f$ where we can specify the parameters $N, \ R, \ S(k)$ in an explicit and uniform way. We concentrate on this problem below.
Now we recall some well-known features of Taylor domination. Basically, it allows us to compare the behavior of $f(z)$ with
the behavior of the polynomial $P_{N}(z)=\sum_{k=0}^{N}a_{k}z^{k}$. In particular, the number of zeroes of $f$ can be easily
bounded in this way. In one direction the bound is provided by the classical result of \cite{biernacki1936fonctions}.
To formulate it, we need the following definition (see \cite{hayman1994multivalent} and references therein):

\begin{definition}
\label{def:pvalent}A function $f$ regular in a domain $\O\subset{\mathbb{C}}$
is called $p$-valent there, if for any $c\in{\mathbb{C}}$ the number
of solutions in $\O$ of the equation $f(z)=c$ does not exceed $p$.
\end{definition}

\begin{theorem}[Biernacki, 1936, \cite{biernacki1936fonctions}]
\label{thm:bier}If $f$ is $p$-valent in the disk $D_{R}$ of radius $R$ centered at $0\in{\mathbb{C}}$ then for each
$k\geq p+1$
\[
|a_{k}|R^{k}\le(A(p)k/p)^{2p}\max_{i=1,\ldots,p}|a_{i}|R^{i},
\]
where $A(p)$ is a constant depending only on $p$.
\end{theorem}
In our notations, \prettyref{thm:bier} claims that a function $f$
which is $p$-valent in $D_{R},$ possesses a $(p,R,(Ak/p)^{2p})$
- Taylor domination property.

For univalent functions, i.e. for $p=1,\ R=1,$ \prettyref{thm:bier} gives $|a_{k}|\le A(1)^2k^2|a_{1}|$
for each $k$, while the sharp bound of the Bieberbach conjecture is $|a_{k}|\le k|a_{1}|$.

Various forms of inverse results to \prettyref{thm:bier} are known. In particular, an explicit bound for the number of zeroes
of $f$ possessing Taylor domination can be obtained by combining \prettyref{prop:Root} and Lemma 2.2 from \cite{roytwarf1997bernstein}:

\begin{theorem}\label{thm:Roy.Yom} Let the function $f$ possess an $(N,R,S(k))$ - Taylor domination property. Then for each
$R'<R$, $f$ has at most $M(N,\frac{R'}{R},S(k))$ zeros in $D_{R'}$, where $M(N,\frac{R'}{R},S(k))$ is a function depending only on $N$,
$\frac{R'}{R}$ and on the sequence $S(k)$, satisfying $\lim_{{{R'}\over R}\to 1}M=\infty$ and $M(N,\frac{R'}{R},S)=N$ for $\frac{R'}{R}$
sufficiently small.
\end{theorem}
We can replace the bound on the number of zeroes of $f$ by the bound on its valency, if we exclude $a_0$ in the definition of
Taylor domination (or, alternatively, if we consider the derivative $f'$ instead of $f$).

\section{Taylor domination and recurrence relations}\label{TD.rec}

We start with a very general result obtained in \cite{roytwarf1997bernstein}.
Assume that a sequence $\Phi$ of mappings $\phi_k:{\mathbb C}^{d+1}\rightarrow {\mathbb C}$
is given, $k=d+1,\ d+2,\cdots$. For any $w=(w_0,\cdots,w_d)\in {\mathbb C}^{d+1}$ construct a sequence $a_k(w)$ as follows: \
$a_i=w_i,\ i=0,\cdots,d$, and $a_j=\phi_j(w)$ for $j>d$. We also consider a power series $f_w(z)=\sum^\infty_{k=0} a_k(w)z^k$.
Of course, any recurrence relation produces such a sequence $\Phi$ by iteration.

\smallskip

Assume also that each $\phi_k$ is a Lipschitzian mapping, satisfying

\begin{equation}
|\phi_k(w)|\leq C^k |w|
\end{equation}
for any $k\geq d+1$ and any $w\in B(0,\delta)$, with some given $\delta>0$ and $C>0$. This is the case in most of natural examples.

\begin{theorem}\label{Roy.Yom.rec}(Theorem 4.1, \cite{roytwarf1997bernstein})
For $C$ and $\delta$ as above and for any $w\in B(0,\delta)$, the series $f_w(x)$ converges on $D_R, \ R=1/C$ and possesses there
$(d,R,K)$-Taylor domination, with $K=[\max(1,C)]^d$.
\end{theorem}
Next we restrict ourselves to linear recurrences, and produce more explicit bounds. We consider the class $\cal S$ of linear
non-stationary homogeneous recurrence relations $\cal R$ of a fixed length, with uniformly bounded coefficients:

\begin{equation}
a_{k}=\sum_{j=1}^{d}c_{j}(k)\cdot a_{k-j},\ \ k=d,d+1,\dots,\label{eq:mainrec}
\end{equation}
If for $j=1,\ldots, d$ the coefficients $c_j(k)$ have a form $c_j(k)=c_{j}+\psi_{j}(k)$, with fixed $c_j$ and with
$\lim_{k\rightarrow\infty}\psi_{j}(k)=0$, then recurrence relation \eqref{eq:mainrec} is said to be
a \emph{linear recurrence relation of Poincaré type }(see \cite{perron1921summengleichungen,poincare1885equations}). We denote the class
of such recurrences $\cal S_P.$

We would like to write the bounds on $c_j(k)$ in a form
\[
\left|c_{j}\left(k\right)\right|\leqslant K\rho^{j},\qquad j=1,\dots,d, \ k=d, d+1, \ldots \ ,
\]
for certain positive constants $K,\rho$. So for each ${\cal R}\in {\cal S}$ we define $K({\cal R})$ and $\rho({\cal R})$ to be the pair
of constants providing the required bounds on $c_j(k)$, for which the product $\nu({\cal R})=(2K({\cal R})+2)\cdot \rho({\cal R})$ is
minimal possible. We put $R({\cal R})\isdef \nu({\cal R})^{-1}.$

\begin{theorem}(Theorem 3.1, \cite{Bat.Yom.2014})\label{thm:main-result}
Let $\left\{ a_{k}\right\} _{k=0}^{\infty}$ be a solution of the recurrence relation
${\cal R}\in {{\cal S}}$. Put $K=K({\cal R}), \ \rho=\rho({\cal R}), \ R=R({\cal R}).$ Then the series $f(z)=\sum_{k=0}^\infty a_kz^k$
converges in the open disk $D_R$ and possesses there
\newline
$(d-1,R,(2K+2)^{d-1})$ Taylor domination.
\end{theorem}
By a proper rescaling, \prettyref{thm:main-result} can be easily extended to non-stationary linear recurrences with a subexponential
(or exponential) growth of the coefficients $c_j(k)$. Consequently, generating functions of such recurrences allow for explicit bounds
on their valency. On the other hand, a drawback of this approach is that in the case of linear recurrences with constant coefficients
(and for Poincar\'e-type recurrences - see below) the disk $D_R$ where the uniform Taylor domination is guaranteed, is much smaller
than the true disk of convergence.

\smallskip

An important feature of \prettyref{thm:main-result} is that it allows us to provide an essentially complete characterisation of
solutions of recurrence relations ${\cal R}\in {\cal S}$ through Taylor domination. The following result is new, although it follows
closely the lines of Theorem 2.3 of \cite{Fri.Yom}. Accordingly, we give only a sketch of the proof, referring the reader to \cite{Fri.Yom}
for details.

\begin{theorem}\label{thm:charact.rec}
A sequence $\left\{ a_{k}\right\} _{k=0}^{\infty}$ is a solution of the recurrence relation ${\cal R}\in {{\cal S}}$ of length $d$ if and
only if for each $m$ its subsequence $\left\{ a_{m+k}\right\} _{k=0}^{\infty}$ possesses a $(d-1,R,C)$ Taylor domination for some positive
$R$ and $C$.
\end{theorem}
\begin{svmultproof2}
In one direction the result follows directly from \prettyref{thm:main-result}. Conversely, if for each $m$ the subsequence
$\left\{ a_{m+k}\right\} _{k=0}^{\infty}$ possesses a $(d-1,R,C)$ Taylor domination, we use the corresponding bound on $|a_{m+d}|$ to
construct step by step the coefficients $c_j(k), \ j=1,\ldots,d$ in ${\cal R},$ in such a way that they remain uniformly bounded in $k$.
\end{svmultproof2}
Let us now consider Poincar\'e-type recurrences. Characterization of their solutions looks a much more challenging problem than the one
settled in \prettyref{thm:charact.rec}. Still, one
can expect deep connections with Taylor domination. One result in this direction (which is a sharpened version of
\prettyref{thm:main-result} above), was obtained in \cite{Bat.Yom.2014}. For ${\cal R}\in {\cal S_P}$ the characteristic polynomial
$\sigma^d-\sum_{j=1}^d c_j\sigma^{d-j}=0$ and the characteristic roots $\sigma_1,\ldots,\sigma_d$ of ${\cal R}$ are those of its
constant part. We put $\rho({\cal R})=\max_j |\sigma_j|.$

\begin{theorem}(Theorem 5.3, \cite{Bat.Yom.2014})\label{thm:Poinc.Rec} Let $\left\{a_{k}\right\} _{k=0}^{\infty}$ satisfy a fixed
recurrence ${\cal R}\in {\cal S_P}$. Put $\rho\isdef \rho({\cal R}), \ R=2^{-(d+3)}\rho^{-1}.$ Let $\hat N$ be the minimal of
the numbers $n$ such that for all $k>n$ we have $|\psi_{j}(k)|\leq 2^d\rho^j, \ j=1,\dots,d$. We put $N=\hat{N}+d,$ and
$C=2^{(d+3)N}.$

Then $\left\{a_{k}\right\} _{k=0}^{\infty}$ possesses $(N,R,C)$-Taylor domination property.
\end{theorem}
For Poincar\'e-type recurrences one can ask for Taylor domination in the maximal disk of convergence, which is typically $\rho^{-1}.$
We discuss this problem in the next section.

\subsection{Tur\'an's lemma, and a possibility of its extension}\label{Turan}

It is well known that the Taylor coefficients $a_k$ of a rational function $R(z)=\frac{P(z)}{Q(z)}$ of degree $d$ satisfy a linear recurrence
relation with constant coefficients
\begin{equation}
a_{k}=\sum_{j=1}^{d}c_{j}a_{k-j},\ \ k=d,d+1,\dots,\label{eq:rec.Const.Coef}
\end{equation}
where $c_j$ are the coefficients of the denominator $Q(z)$ of $R(z)$. Conversely, for any initial terms $a_0,\ldots,a_{d-1}$ the
solution sequence of \eqref{eq:rec.Const.Coef} forms a sequence of the Taylor coefficients $a_k$ of a rational function $R(z)$
as above. Let $\sigma_1,\ldots,\sigma_d$ be the characteristic roots of \eqref{eq:rec.Const.Coef}, i e. the roots of its characteristic
equation $\sigma^d-\sum_{j=1}^d c_j\sigma^{d-j}=0.$

\smallskip

Taylor domination property for rational functions is provided by the following theorem, obtained in \cite{Bat.Yom.2014}, which is,
essentially, equivalent to the ``first Turán lemma'' (\cite{turan1953neue,turan1984new,Ineq}, see also \cite{nazarov1994local}):

\begin{theorem}\label{thm:turan}(Theorem 3.1, \cite{Bat.Yom.2014}) Let $\left\{ a_{j}\right\} _{j=1}^{\infty}$ satisfy
recurrence relation \eqref{eq:rec.Const.Coef} and let $\sigma_1,\ldots,\sigma_d$ be its characteristic roots. Put
$R\isdef\min_{i=1,\dots d}\left|\sigma_{i}^{-1}\right|.$ Then for each $k\geq d$

\begin{equation}
\left|a_{k}\right|R^{k}\leq \ Q(k,d)\max_{i=0,\dots,d-1}\ |a_{i}|R^{i},\label{eq:tur}
\end{equation}
where $Q(k,d)=[2e({k\over d}+1)]^d$.
\end{theorem}
\prettyref{thm:turan} provides a uniform Taylor domination for rational functions in their maximal disk of convergence $D_{R}$,
in the strongest possible sense. Indeed, after rescaling to the unit disk $D_{1}$ the parameters of \eqref{eq:tur} depend only on
the degree of the function, but not on its specific coefficients.

\smallskip

We consider a direct connection of Turán's lemma to Taylor domination, provided by \prettyref{thm:turan} as an important and promising
fact. Indeed, Turán's lemma has numerous applications and conncetions, many of them provided already in the first Turan's book
\cite{turan1953neue}. It can be considered as a result on exponential polynomials, and in this form it was a starting point for many deep
investigations in Harmonic Analysis, Uncertainty Principle, Analytic continuation, Number Theory
(see \cite{Ineq,nazarov1994local,turan1953neue,turan1984new} and references therein). Recently some applications in Algebraic Sampling
were obtained, in particular, estimates of  robustness of non-uniform sampling of ``spike-train'' signals
(\cite{Sampl,friedland2011observation}). One can hope that apparently new connections of Turán's lemma with Taylor domination,
presented in \cite{Bat.Yom.2014} and in the present paper, can be further developed.

\smallskip

A natural open problem, motivated by \prettyref{thm:turan}, is a possibility to extend uniform Taylor domination in the maximal disk
of convergence $D_{R}$, as provided by \prettyref{thm:turan} for rational functions, to wider classes of generating functions of
Poincaré type recurrence relations. Indeed, for such functions the radius of convergence of the Taylor series is, essentially, the same
as for the constant-coefficients recurrences - it is the inverse of one of the characteristic roots: for some $\sigma_{j}$

\[
\limsup_{k\to\infty}\sqrt[k]{\left|a_{k}\right|}=\left|\sigma_{j}\right|.
\]
O.Perron proved in \cite{perron1921summengleichungen} that this relation
holds for a general recurrence of Poincaré type, but with an additional
condition that $c_{d}+\psi_{d}\left(k\right)\neq0$ for all $k\in\N$.
In \cite{pituk2002more} M.Pituk removed this restriction, and proved
the following result.

\begin{theorem}[Pituk's extension of Perron's Second Theorem, \cite{pituk2002more}]
\label{thm:pituk-scalar}Let $\left\{ a_{k}\right\} _{k=0}^{\infty}$
be any solution to a recurrence relation ${\cal R}$ of Poincaré class
$\SP$. Then either $a_{k}=0$ for $k\gg1$ or
\[
\limsup_{k\to\infty}\sqrt[k]{\left|a_{k}\right|}=\left|\sigma_{j}\right|,
\]
where $\sigma_{j}$ is one of the characteristic roots of ${\cal R}$.
\end{theorem}
This result implies the following:
\begin{theorem}(Theorem 5.2, \cite{Bat.Yom.2014})\label{thm:poinc-dom-maximal-disk}
Let $\left\{ a_{k}\right\} _{k=0}^{\infty}$
be any nonzero solution to a recurrence relation ${\cal R}$ of Poincaré class ${\cal S_P}$ with initial data $\bar{a}$, and let
$R$ be the radius of convergence of the generating function $f\left(z\right)$. Then necessarily $R>0$, and in fact
$R=\left|\sigma\right|^{-1}$ where $\sigma$ is some (depending on $\bar{a}$) characteristic root of ${\cal R}$. Consequently, $f$
satisfies $\left(d-1,R,S\left(k\right)\right)$-Taylor domination with $S\left(k\right)$ as defined in \prettyref{prop:Root}.
\end{theorem}
Taylor domination in the maximal disk of convergence provided by \prettyref{thm:poinc-dom-maximal-disk}, is only partially effective.
Indeed, the number $d-1$ and the radius $R$ are as prescribed by the constant part of the recurrence. However, \prettyref{prop:Root}
only guarantees that the sequence $S\left(k\right)=R^k|a_k|\cdot (\max_{i=0}^{d-1}R^i|a_i|)^{-1}$ is of
subexponential growth but gives no further information on it. We can pose a natural question in this direction. For a sequence
$\Delta=\{\delta_k\}$ of positive numbers tending to zero, consider a subclass $\SPD$ of $\SP$,
consisting of ${\cal R}\in \SP$ with $|\psi_{j}(k)|\leq \delta_k\cdot \rho({\cal R})^j, \ j=1,\ldots,d, \ k=d,d+1,\ldots$

\smallskip

\noindent{\bf Problem 1.} {\it Do solutions of recurrence relations ${\cal R}\in \SPD$ possess
$(N,R,S(k))$-Taylor domination in the maximal disk of convergence $D_R$, with $S(k)$ depending only on $d$ and $\Delta$?
Is this true for specific $\Delta$, in particular, for $\Delta= \{1,{1\over 2},{1\over 3},...,\},$ as it occurs in most of
examples (solutions of linear ODE's, etc.)?}

\medskip

Taking into account well known difficulties in the analysis of Poincaré-type recurrences, this question may be tricky.
Presumably, it can be easier for $\Delta$ with $\sum_{k=1}^\infty \delta_k < \infty$.

\smallskip

Some initial examples in this direction were provided in \cite{yomdin2014Bautin}, via techniques of Bautin ideals. In Section
\ref{Bautin} below we indicate an interrelation of the Poincaré-type recurrences with the Bautin ideals techniques. Another possible
approach may be via an inequality for consecutive moments of linear combinations of $\delta$-functions, provided by Theorem 3.3 and
Corollary  3.4 of \cite{yomdin2010sing.Prony}. This inequality is closely related to Turán's lemma (and to the ``Tur\'an third theorem''
of \cite{turan1953neue,turan1984new,Ineq}). It was obtained via techniques of finite differences which, presumably, can be extended to
Stieltjes transforms of much wider natural classes of functions, like piecewise-algebraic ones.

One can consider some other possible approaches to the ``Tur\'an-like'' extension of Taylor domination to the full
disk of convergence $D_{R}$ for the Poincaré-type recurrences. First, asymptotic expressions in \cite{bodine2004asymptotic,pituk2002more}
may be accurate enough to provide an inequality of the desired form. If this is the case, it remains to get explicit bounds in
these asymptotic expressions.

Second, one can use a ``dynamical approach'' to recurrence relation
\eqref{eq:mainrec} (see \cite{borcea2011parametric,coppel1971dichotomies,kloeden2011non,potzsche2010geometric,yom.nonaut.dyn}
and references therein). We consider \eqref{eq:mainrec} as a non-autonomous
linear dynamical system $T$. A ``non-autonomous diagonalization'' of $T$ is a sequence ${\cal L}$ of linear changes of variables,
bringing this system to its ``constant model'' $T_{0}$, provided by the limit recurrence relation \eqref{eq:rec.Const.Coef}.

If we could obtain a non-autonomous diagonalization ${\cal L}$ of $T$ with an explicit bound on the size of the linear changes of
variables in it, we could get the desired inequality as a pull-back, via ${\cal L}$, of the Turán inequality for $T_{0}$. There
are indications that this approach may work in the classes $\SPD$ with $\sum_{k=1}^\infty \delta_k < \infty$.

\smallskip

The following partial result in the direction of Problem 1 above was obtained in \cite{Bat.Yom.2014}. It provides Taylor domination
{\it in a smaller disk}, but with explicit parameters, expressed in a transparent way through the constant part of ${\cal R}$, and
through the size of the perturbations.

\begin{theorem}(Corollary 5.1, \cite{Bat.Yom.2014}) \label{thm:Poinc.Rec.Cor}
Let $\Delta=\{\delta_k\}$ be a sequence of positive numbers tending to zero. Define $\hat N(\Delta)$ as a minimal number $n$
such that for $k>n$ we have $\delta_k\leq 2^d$. Then for each ${\cal R}\in \SPD$, the solution sequences of ${\cal R}$
possess $(N,R,C)$-Taylor domination, where $N=\hat N(\Delta)+d, \ R=2^{-(d+3)}\rho({\cal R})^{-1}, \ C=2^{(d+3)N}.$
\end{theorem}

\subsection{An example: D-finite functions}

\global\long\def\np{p}
\global\long\def\ff#1#2{\left(#1\right)_{#2}}
\global\long\def\mvi{\ensuremath{\Lambda}}

In this section we briefly summarize results of \cite{Bat.Yom.2014}, concerning a certain class of power series, defined by
the Stieltjes integral transforms
\begin{equation}
f\left(z\right)=S_{g}\left(z\right)=\int_a^b\frac{g\left(x\right)\D x}{1-zx},\label{eq:def-stieltjes}
\end{equation}
where $g\left(x\right)$ belongs to the class ${\cal PD}$ of the so-called \emph{piecewise
D-finite functions} \cite{bat2008}, which are solutions of linear ODEs with polynomial coefficients, possessing a finite number of
discontinuities of the first kind.

Using the expansion $\left(1-zx\right)^{-1}=\sum_{k=0}^{\infty}\left(zx\right)^{k}$ for $\left|z\right|<\frac{1}{\left|x\right|}$,
we obtain the following useful representation of $S_{g}\left(z\right)$ through the moments of $g$:
\[
S_{g}\left(z\right)=\sum_{k=0}^{\infty}m_{k}z^{k},\quad\text{where }m_{k}\isdef\int_{a}^{b}x^{k}g\left(x\right)\D x.
\]
Obtaining uniform Taylor domination for $S_{g},$ where $g$ belongs to particular subclasses of ${\cal PD}$ (in particular, $g$
being piecewise algebraic), is an important problem with direct applications in Qualitative Theory of ODEs
(see \cite{Bli.Bri.Yom,briskin2010center} and references therein).

\begin{definition}
A real-valued bounded integrable function $g:\left[a,b\right]\to\R$
is said to belong to the class ${\cal PD}\left(\Op,\np\right)$ if
it has $0\leqslant\np<\infty$ discontinuities (not including the
endpoints $a,b$) of the first kind, and between the discontinuities
it satisfies a linear homogeneous ODE with polynomial coefficients
$\Op g=0$, where
\[
\Op=\sum_{j=0}^{n}p_{j}\left(x\right)\left(\frac{\D}{\D x}\right)^{j},\quad p_{j}\left(x\right)=\sum_{i=0}^{d_{j}}a_{i,j}x^{i}.
\]

\end{definition}
Let $g\in{\cal PD}\left(\Op,\np\right)$, with $\Op$ as above. Denote the discontinuities of $g$ by
$a=x_{0}<x_{1}<\dots<x_{\np}<x_{\np+1}=b$. In what follows, we shall use some additional notation. Denote for each $j=0,\dots,n,$
\ \ $\alpha_{j}\isdef d_{j}-j.$ Let $\alpha\isdef\max_{j}\alpha_{j}$.

Our approach in \cite{Bat.Yom.2014} is based on the following result:
\begin{theorem}[\cite{bat2008}] Let $g\in{\cal PD}\left(\Op,\np\right)$. Then the moments $m_{k}=\int_{a}^{b}g\left(x\right)\D x$
satisfy the recurrence relation

\begin{equation}
\sum_{\ell=-n}^{\alpha}q_{\ell}\left(k\right)m_{k+\ell}=\e_{k},\quad k=0,1,\dots,\label{eq:moments-main-rec}
\end{equation}
where $q_{\ell}\left(k\right)$ are polynomials in $k$ expressed through the coefficients of $\Op$, while $\e_k$ in the right hand
side are expressed through the values of the coefficients of $\Op$, and through the jumps of $g$, at the points
$x_{0}<x_{1}<\dots<x_{\np}<x_{\np+1}$.
\end{theorem}
The recurrence \eqref{eq:moments-main-rec} is inhomogeneous, and the coefficient of the highest moment may vanish for some $k$.
Accordingly, we first transform \eqref{eq:moments-main-rec} into a homogeneous matrix recurrence.

\begin{definition}
The vector function $\vec{y}\left(k\right):\N\to\C^{n}$ is said to
satisfy a linear system of Poincaré type, if
\begin{equation}
\vec{y}\left(k+1\right)=\left(A+B\left(k\right)\right)\vec{y}\left(k\right),\label{eq:system-poincare}
\end{equation}
where $A$ is a constant $n\times n$ matrix and $B\left(k\right):\N\to\C^{n\times n}$
is a matrix function satisfying $\lim_{k\to\infty}\|B\left(k\right)\|=0$.
\end{definition}
Put $\tau\isdef n\left(\np+2\right)$. Now define the vector function $\vec{w}\left(k\right):\N\to\C^{\alpha+n+\tau}$ as
\[
\vec{w}\left(k\right)\isdef\begin{bmatrix}m_{k-n}\\
\vdots\\
m_{k+\alpha-1}\\
\e_{k}\\
\vdots\\
\e_{k+\tau-1}
\end{bmatrix}.
\]
\begin{proposition}\label{prop:System}
The vector function $\vec{w}\left(k\right)$ satisfies a linear system of the form \eqref{eq:system-poincare}. This system
is of Poincaré type, if and only if
\begin{equation}
\alpha_{n}\geqslant\alpha_{j}\qquad j=0,1,\dots,n.\label{eq:poinc-sys-cond}.
\end{equation}
This last condition is equivalent to the operator $\Op$ having at most a regular singularity at $z=\infty$. The set $Z_{A}$ of
the eigenvalues of the matrix $A$ is precisely the union of the roots of $p_{n}\left(x\right)$ (i.e. the singular points of
the operator $\Op$) and the jump points $\left\{ x_{i}\right\} _{i=0}^{\np+1}$.
\end{proposition}
Now we establish in \cite{Bat.Yom.2014} Taylor domination for the Stieltjes transform $S_{g}\left(z\right)$, combining two
additional results: the first is the system version of \prettyref{thm:pituk-scalar}.

\begin{theorem}[\cite{pituk2002more}]
\label{thm:pituk-system}Let the vector $\vec{y}\left(k\right)$
satisfy the perturbed linear system of Poincaré type \eqref{eq:system-poincare}.
Then either $\vec{y}\left(k\right)=\vec{0}\in\C^{n}$ for $k\gg1$
or
\[
\lim_{k\to\infty}\sqrt[k]{\|\vec{y}\left(k\right)\|}
\]
exists and is equal to the modulus of one of the eigenvalues of the
matrix $A$.
\end{theorem}
Next main problem is: \emph{how many first moments $\left\{ m_{k}\right\} $ can vanish for a nonzero
$g\in{\cal PD}\left(\Op,\np\right)$?}. $N$ in Taylor domination cannot be smaller than this number.
In \cite{batbinZeros} we study this question, proving the following result.
\begin{theorem}[\cite{batbinZeros}]
\label{thm:moment-vanishing}Let the operator $\Op$ be of Fuchsian
type (i.e. having only regular singular points, possibly including
$\infty$). In particular, $\Op$ satisfies the condition \eqref{eq:poinc-sys-cond}.
Let $g\in{\cal PD}\left(\Op,\np\right)$.
\begin{enumerate}
\item If there is at least one discontinuity point $\xi\in\left[a,b\right]$
of $g$ at which the operator $\Op$ is nonsingular, i.e. $p_{n}\left(\xi\right)\neq0$,
then vanishing of the first $\tau+d_{n}-n$ moments $\left\{ m_{k}\right\} _{k=0}^{\tau+d_{n}-n-1}$
of $g$ implies $g\equiv0$.
\item Otherwise, let $\Lambda\left(\Op\right)$ denote the largest positive
integer characteristic exponent of $\Op$ at the point $\infty$.
In fact, the indicial equation of $\Op$ at $\infty$ is $q_{\alpha}\left(k\right)=0$.
Then the vanishing of the first $\Lambda\left(\Op\right)+1+d_{n}-n$
moments of $g$ implies $g\equiv0$.
\end{enumerate}
\end{theorem}
With these theorems in place, the following result is obtained in \cite{Bat.Yom.2014}:
\begin{theorem}
\label{thm:moment-taylor-dom}Let $g\in{\cal PD}\left(\Op,\np\right)$
be a not identically zero function, with $\Op$ of Fuchsian type.
Then the Stieltjes transform $S_{g}\left(z\right)$ is analytic at
the origin, and the series
\[
S_{g}\left(z\right)=\sum_{k=0}^{\infty}m_{k}z^{k}
\]
converges in a disk of radius $R$ which satisfies
\[
R\geqslant R^{*}\isdef\min\left\{ \xi^{-1}:\;\xi\in Z_{A}\right\} ,
\]
where $Z_{A}$ is given by \prettyref{prop:System}. Furthermore, for every
\[
N\geqslant\max\left\{ \tau-1,\Lambda\left(\Op\right)\right\} +d_{n}-n,
\]
$S_{g}$ satisfies $\left(N,R,S\left(k\right)\right)$ Taylor domination,
where $S\left(k\right)$ is given by \prettyref{prop:Root}. \end{theorem}

\section{Bautin Ideals}\label{Bautin}

In this section we consider families

\begin{equation}\label{eq:bautin}
f_\lambda(z)=\sum_{k=0}^{\infty}a_{k}(\lambda)z^{k},
\end{equation}
where $\lambda = (\lambda_1,\ldots,\lambda_n)\in {\mathbb C}^n$. We assume that each $a_{k}(\lambda)$ is a polynomial in $\lambda$.
It was a remarkable discovery of N. Bautin (\cite{Bau1,Bau2}) that in this situation the behavior of the Taylor coefficients
of $f_\lambda$, and consequently, of its zeroes, can be understood in terms of the ideals
$I_k = \{a_{0}(\lambda),\ldots,a_{k}(\lambda)\}$ generated by the subsequent Taylor coefficients $a_j(\lambda)$ in the
polynomial ring ${\mathbb C}[\lambda]$. Explicit computation of $I_k$ in specific examples may be very difficult. In particular,
Bautin himself computed $I_k$ for the Poincar\'e first return mapping of the plane polynomial vector field of degree $2$,
producing in this way one of the strongest achievements in the Hilbert 16th problem up to this day: at most three limit cycles can
bifurcate from an isolated equilibrium point in vector field of degree $2$.

\smallskip

Bautin's approach can be extended to a wide classes of families of the form \eqref{eq:bautin}. In particular, the following class
was initially defined and investigated in \cite{Bri.Yom,Fra.Yom,roytwarf1997bernstein}, and further studied in
\cite{Yom.Baut,yomdin2014Bautin}:

\begin{definition}\label{A0.series}
Let $f_\lambda(z)=\sum_{k=0}^{\infty}a_{k}(\lambda)z^{k}$. \ $f$ is called an $A_0$-series if the following condition is satisfied:

$$
\deg a_k(\lambda) \leq K_1k+K_2, \ \ \text{and} \ \ ||a_k(\lambda)||\leq K_3\cdot K_4^k, \ k=0,1,\ldots,
$$
for some positive $K_1,K_2,K_3,K_4$.

\smallskip

The ideal $I=\{a_{0}(\lambda),\ldots,a_{k}(\lambda),\ldots,\}$ generated by all the subsequent Taylor coefficients $a_j(\lambda)$
in the polynomial ring ${\mathbb C}[\lambda]$ is called the Bautin ideal of $f$, and the minimal $d$ such that
$I=I_d=\{a_{0}(\lambda),\ldots,a_{d}(\lambda)\}$ is called the Bautin index of $f$.
\end{definition}
Such a finite $d$ exists by the Noetherian property of the ring ${\mathbb C}[\lambda]$.

The following result of \cite{Fra.Yom} connects the algebra of the Bautin ideal of $f$ with Taylor domination for this series:

\begin{theorem}\label{thm:Fra.Yom}
Let $f_\lambda(z)=\sum_{k=0}^{\infty}a_{k}(\lambda)z^{k}$ be an $A_0$-series, and let $I$ and $d$ be the Bautin ideal and the
Bautin index of $f$. Then $f$ possesses a $(d,R,C)$-Taylor domination for some positive $R,C$ depending on $K_1,K_2,K_3,K_4$,
and on the basis $a_{0}(\lambda),\ldots,a_{d}(\lambda)$ of the ideal $I$.
\end{theorem}
\begin{svmultproof2}
The idea of the proof is very simple: if for the $A_0$-series $f_\lambda(z)=\sum_{k=0}^{\infty}a_{k}(\lambda)z^{k}$ the first
$d$ Taylor coefficients $a_{0}(\lambda),\ldots,a_{d}(\lambda)$ generate the ideal $I=\{a_{0}(\lambda),\ldots,a_{k}(\lambda),\ldots,\}$,
then for any $k>d$ we have

\[
a_k(\lambda)=\sum_{i=0}^d \psi_i^k(\lambda)a_i(\lambda),
\]
with $\psi_i^k(\lambda)$ - certain polynomials in $\lambda$. The classical Hironaka's division theorem provides a bound on the
degree and the size of $\psi_i^k(\lambda)$, and hence we get an explicit bound on $|a_k(\lambda)|$ through
$\max \{|a_{0}(\lambda)|,\ldots,|a_{d}(\lambda)|\}$ for each $\lambda$.
\end{svmultproof2}
Also here to get an explicit and uniform in $\lambda$ Taylor domination in the full disk of convergence is a difficult problem.
Some general additional conditions on $A_0$-series, providing such result, were given in \cite{Yom.Baut}. These conditions are
technically rather involved. A simpler special case was treated in \cite{yomdin2014Bautin}.

\subsection{Recurrence relations and Bautin ideals}

A rather detailed investigation of $A_0$-series produced by recurrence relations was provided in \cite{Bri.Yom}. It was done in
several situations, including algebraic functions, linear and nonlinear differential equations with polynomial coefficients, and
the Poincar\'e first return mapping for the Abel differential equation. In all these cases, except the last one, it was shown
that the resulting series are $A_0$-ones, and their Bautin ideal was computed. In the present paper we generalize the result of
\cite{Bri.Yom}, showing how to compute the Bautin ideal for general non-stationary polynomial recurrences. Let, as above,
$\lambda = (\lambda_1,\ldots,\lambda_n)\in {\mathbb C}^n$, and put $u=(u_1,\ldots,u_d)\in {\mathbb C}^d.$ Consider a sequence
of polynomials
$$
P_k(u_1,\ldots,u_d)=\sum_{|\alpha|\leq d_k} A_{k,\alpha}(\lambda)u^\alpha,
$$
where
$$
A_{k,\alpha}(\lambda)=\sum_{|\beta|\leq d_{k,\alpha}} A_{k,\alpha,\beta}\lambda^\beta,
$$
are, in their turn, polynomials in $\lambda$ of degrees $d_{k,\alpha}$. Here
$\alpha=(\alpha_1,\ldots,\alpha_d)\in {\mathbb N}^d$, and $\beta=(\beta_1,\ldots,\beta_n)\in {\mathbb N}^n$ are
multi-indices, $|\alpha|=|\alpha_1|+\ldots+|\alpha_d|,$ \ $|\beta|=|\beta_1|+\ldots+|\beta_n|,$ \
$u^\alpha=u_1^{\alpha_1}\cdot\ldots \cdot u_d^{\alpha_d},$ and
$\lambda^\beta=\lambda_1^{\beta_1}\cdot\ldots \cdot \lambda_n^{\beta_n}.$ Finally, $A_{k,\alpha,\beta}$ are complex constants.

\begin{definition}\label{Non.Stat.Polyn.Rec}
A recurrence of the form
\begin{equation}\label{eq:Polyn.Rec}
a_k(\lambda)=P_k(a_{k-1}(\lambda),\ldots,a_{k-d}(\lambda)), \ k=d,d+1,\dots
\end{equation}
is called a non-stationary polynomial recurrence relation of length $d$.
\end{definition}
The following result is purely algebraic, so we put no restrictions on the degrees and the size of the polynomials $P_k$.

\begin{theorem}\label{thm:Baut.Rec.General}
Let polynomials $a_0(\lambda),\ldots,a_{d-1}(\lambda)$ in $\lambda$ be given, and let $a_k(\lambda) , \ k=d,d+1,\dots$ be
produced by recurrence \eqref{eq:Polyn.Rec}. Then the Bautin ideal $I$ of the formal series
$f_\lambda(z)=\sum_{k=0}^\infty a_k(\lambda)z^k$ is generated by $a_0(\lambda),\ldots,a_{d-1}(\lambda).$
\end{theorem}
\begin{svmultproof2}
Equation \eqref{eq:Polyn.Rec} applied with $k=d$ shows that $a_d(\lambda)$ belongs to the ideal $I_{d-1}$ generated by
$a_0(\lambda),\ldots,a_{d-1}(\lambda).$ Applying this equation step by step, we see that all $a_k(\lambda), \ k=d,d+1,\ldots$
belong to $I_{d-1}$. Therefore $I=I_{d-1}$.
\end{svmultproof2}
We expect that if we assume the degrees $d_k$ and $d_{k,\alpha}$, as well as the coefficients $A_{k,\alpha,\beta}$ to be
uniformly bounded, then $f_\lambda(z)$ is, in fact, an $A_0$-series. We expect that this fact can be shown by the methods
of \cite{Bri.Yom}, properly extended to non-stationary polynomial recurrence relation. Let us consider an example in this
direction.

\smallskip

Let
$$
P_k(u_1,\ldots,u_d)=\sum^d_{j=1} A_{k,j}(\lambda)u_j, \ \ A_{k,j}(\lambda)=\sum_{i=1}^n A_{k,j,i}\lambda^i, \ \
A_{k,j,i}\in {\mathbb C}.
$$
We consider linear recurrence

\begin{equation}\label{eq:lin.lambda}
a_k(\lambda)=\sum^d_{j=1} A_{k,j}(\lambda)a_{k-j}(\lambda).
\end{equation}

\begin{proposition}\label{prop:Baut.Rec.General}
Let polynomials $a_0(\lambda),\ldots,a_{d-1}(\lambda)$ of degrees $0,1,\ldots,d-1$ in $\lambda$ be given, and let
$a_k(\lambda) , \ k=d,d+1,\dots$ be produced by recurrence \eqref{eq:lin.lambda}, with $A_{k,j,i}$ uniformly bounded
in $k$. Then the formal series $f_\lambda(z)=\sum_{k=0}^\infty a_k(\lambda)z^k$ is, in fact, an $A_0$ series. In particular,
the degree of $a_k(\lambda)$ is at most $k$, and the Bautin ideal $I$ is generated by $a_0(\lambda),\ldots,a_{d-1}(\lambda).$
\end{proposition}
\begin{svmultproof2}
Equation \eqref{eq:lin.lambda} shows that the degree in $\lambda$ of $a_k$ is by one higher than the maximum of the degrees
of $a_{k-1},\ldots,a_{k-d}$. Applying induction we see that the degree of $a_k(\lambda)$ is at most $k$. The last statement
of \prettyref{prop:Baut.Rec.General} follows directly from \prettyref{thm:Baut.Rec.General}. Now, write

$$
a_k(\lambda)=\sum_{|\beta|\leq k}a_{k,\alpha}\lambda^\beta.
$$
Then equation \eqref{eq:lin.lambda} shows that $a_{k,\beta}$ satisfy a matrix recurrence relation, which in coordinates
takes a form

\begin{equation}\label{eq:coeff.ak}
a_{k,\beta}=\sum_{j=1}^d \sum_{i=1}^n A_{k,j,i} a_{k-j,\beta [i]}.
\end{equation}
Here for $\beta=(\beta_1,\ldots,\beta_n)$ we define $\beta[i]$ as
$$
\beta [i]=(\beta_1,\ldots, \beta_{i-1},\beta_i-1,\beta_{i+1},\ldots,\beta_n).
$$
So in the right hand side of \eqref{eq:coeff.ak} appear all the $a_{k',\beta'}$ with $k'$ between $k-1$ and $k-d$, and
$\beta'$ smaller than $\beta$ by one in exactly one coordinate.

\smallskip

Applying to the matrix recurrence \eqref{eq:coeff.ak} straightforward estimates, we obtain exponential in $k$ upper bounds
on the coefficients $a_{k,\beta}$, so $f_\lambda(z)=\sum_{k=0}^\infty a_k(\lambda)z^k$ is indeed an $A_0$ series.
\end{svmultproof2}
We expect that the Taylor domination bounds of Section \ref{TD.rec} above, as applied to matrix recurrence \eqref{eq:coeff.ak}
can be combined with algebraic consideration of the Bautin ideals. In this way one can hope to provide much more accurate
parameters in Taylor domination than the general ones of \prettyref{thm:Fra.Yom}.

\subsection{Poincar\'e coefficients of Abel equation}\label{Abel}

Consider Abel differential equation

\begin{equation} \label{Abel.de}
y'=p(x)y^2 + q(x) y^3
\end{equation}
with polynomial coefficients $p,q$ on the interval $[a,b]$. A solution $y(x)$ of (\ref{Abel.de}) is called ``closed'' if
$y(a)=y(b)$. The Smale-Pugh problem, which is a version of the (second part of) Hilbert's 16-th problem, is to bound the number
of isolated closed solutions of (\ref{Abel.de}) in terms of the degrees of $p$ and $q$.

\smallskip

This problem can be naturally expressed in terms of the Poincar\'e ``first return'' mapping $y_b=G(y_a)$ along $[a,b]$. Let
$y(x,y_a)$ denote the solution $y(x)$ of (\ref{Abel.de}) satisfying $y(a)=y_a$. The Poincar\'e mapping $G$ associates to each
initial value $y_a$ at $a$ the value $y_b$ at $b$ of the solution $y(x,y_a)$ analytically continued along $[a,b]$. Closed
solutions correspond to the fixed points of $G$. So the problem is reduced to bounding the number of the fixed points of $G$, or
of zeroes of $G(y)-y$. Historically, one of the most successful directions in the study of the Poincar\'e mapping $G$ was the
direction initiated by Bautin in \cite{Bau1,Bau2}: to derive the analytic properties of $G$, in particular, the number of its
fixed points, from the structure of its Taylor coefficients.

\smallskip

It is well known that $G(y)$ for small $y$ is given by a convergent power series

\begin{equation}\label{Poic.map}
G (y)= y + \sum_{k=2}^\infty v_k(p,q,a,b)y^k.
\end{equation}
The ``Poincar\'e coefficients'' $v_k(p,q,x,a)$ of the Poincar\'e mapping from $x$ to $a$ satisfy the following
differential-recurrence relation (see, for example, \cite{briskin2010center}):

\begin{equation}\label{eq:Poic.rec}
{{dv_k}\over {dx}}= -(k-1)pv_{k-1}-(k-2)qv_{k-2}, \ v_0\equiv 0, \ \ v_1\equiv 1,
\ v_k(0)=0, \ k\geq 2.
\end{equation}
This recurrence is apparently not of the form considered above, i.e. it is not ``polynomial recurrence''. Still, it is easy
to see from \eqref{eq:Poic.rec} that the Poincar\'e coefficients $v_k(p,q,b,a)$ are polynomials with rational coefficients in the
parameters of the problem (i.e. in the coefficients of $p$ and $q$). They can be explicitly computed for as large indices $k$ as
necessary. However, because of the derivative in the left hand side, \eqref{eq:Poic.rec} does not preserve ideals. So
\prettyref{thm:Baut.Rec.General} is not applicable to it, and computing Bautin ideals for the Poincar\'e mapping $G(y)$ is in
general a very difficult problem. Also Taylor domination for $G(y)$ is not well understood, besides some special examples. The
only general result concerning the Bautin ideals of $G(y)$ we are aware of was obtained in \cite{Bli.Bri.Yom} using an
approximation of the ``Poincar\'e coefficients'' $v_k(p,q,b,a)$ with certain moment-like expressions of the form
$m_k(p,q,b,a)=\int_a^b P^k(x)q(x)dx, \ P=\int p$. We would like to pose the investigation of the recurrence relation
\eqref{eq:Poic.rec} in the lines of the present paper as an important open problem.

\bibliographystyle{plain}

\end{document}